\documentclass[12pt]{amsart}

\usepackage[margin=1.5in]{geometry}

\usepackage{amssymb}

\usepackage{amsthm}
\usepackage{amsmath}
\usepackage{amsaddr}

\newtheorem{theorem}{Theorem}[section]
\newtheorem{lemma}[theorem]{Lemma}

\theoremstyle{definition}

\newtheorem{corollary}[theorem]{Corollary}

\theoremstyle{remark}

\newcommand{\R}{ {\mathbb R} }

\newcommand{\Z}{ {\mathbb Z} }
\newcommand{\Q}{ {\mathbb Q} }
\newcommand{\gen}[1]{\left\langle #1 \right\rangle}
\newcommand{\ind}[2] {\left( \mathbf{1}_{#1} \right)^{#2}}
\newcommand{\ip}[2]{\left\langle #1,#2 \right\rangle}

\begin{document}

\title[Representations on Simplicial Fans]{Representations of cyclic groups acting on complete simplicial fans}

\author{Jonathan Browder}
\email{browder@math.wustl.edu}
\address{Department of Mathematics, Washington Univeristy in St Louis, Campus Box 1146, St~Louis, MO 63130}
\date{30 October 2008}

\begin{abstract}Let $\sigma$ be a complete simplicial fan in finite dimensional real Euclidean space $V$, and let $G$ be a cyclic subgroup of $GL(V)$ which acts properly on $\sigma$. We show that the representation of $G$ carried by the cohomology of $X_{\sigma}$, the toric variety associated to $\sigma$, is a permutation representation.
\end{abstract}

\maketitle

\section{Introduction}
	Let $L$ be a lattice in $d$-dimensional real Euclidean space $V$, and let $\sigma$ be a polyhedral decomposition of $V$. If the cones of $\sigma$ are generated by vectors in $L$, $\sigma$ is called a (complete) \emph{fan} in $L$ (and we say $\sigma$ is \emph{integral} with respect to $L$). If each $i$-dimensional cone of $\sigma$ is generated by $i$ vectors, we say $\sigma$ is \emph{simplicial}. Associated to any fan $\sigma$ in $V$ is a \emph{toric variety}, $X_{\sigma}$.
  
	In \cite{MR1279220}, Stembridge showed that when $\sigma_R$ is the simplicial fan defined by the hyperplanes orthogonal to the roots of a crystallographic root system $R$, and $G$ is the Weyl group of $R$, the representation of $G$ carried by $H^*(X_{\sigma_R}, {\mathbb Q} )$ is isomorphic to a permutation representation of $G$ (ignoring the grading in the case that $G$ is not of type A or B). The main purpose of this paper is to show an analogous result for cyclic groups $G$ acting appropriately on arbitrary fans.

	Suppose $G$ is a finite subgroup of $GL(V)$ that acts as a group of homomorphisms of $L$ and additionally as a set of automorphisms of some complete simplicial fan $\sigma$. Then $G$ acts on the toric variety $X_{\sigma}$. We say that the action of $G$ is \emph{proper} if whenever an element of $G$ fixes a cone in $\sigma$, then it fixes that cone pointwise.

\begin{theorem} \label{main} Let $G$ be a finite cyclic subgroup of $GL(V)$ acting as a set of homomorphisms of lattice $L$ in $V$ and acting properly on some complete simplicial fan $\sigma$ (integral with respect to $L$). Then the representation of $G$ carried by $H^*(X_{\sigma}, {\mathbb Q} )$ is isomorphic to a permutation representation.
\end{theorem}

(Note that we vary our assumptions slightly from those in section 1 of \cite{MR1279220}, replacing the requirement that $G$ be a subgroup of the orthogonal group with our assumption that $G$ preserves the lattice. Both, of course, hold in the case of the Weyl groups and associated weight lattices under consideration in that paper.)

	For $\sigma$ a complete simplicial fan, let $I$ be the set of primitive vectors of $L$ lying in the one dimensional cones of $\sigma$. We then have a simplicial complex, $\Delta$, associated to $\sigma$ whose faces are the subsets of $I$ whose elements generate a cone in $\sigma$ (note that this is combinatorially equivalent to the complex obtained by intersecting the fan with a sphere centered at the origin in $V$).

Following Stembridge we make use of the fact, due to Danilov \cite{dan:torvar}, that  $H^*(X_{\sigma}, {\mathbb Q} )$ is isomorphic as a graded $\Q G$-module  to $\Q [\Delta]/ \Theta$ (after a degree shift), where $\Q [\Delta]$ is the Stanley-Reisner ring of $\Delta$ and $\Theta$ is a certain system of parameters. As in \cite{MR1279220} we define the \emph{graded character} of a graded $KG$-module $M = \bigoplus_{i} M_i$ to be $\chi [M,q](g) = \sum_i \text{tr}_{M_i}(g)q^i$. If $M$ is finite dimensional as a K-vector space we also have the \emph{ungraded} character of $M$ (i.e., the character of $M$ with the grading disregarded) which is simply $\chi [M,1](g)$. Thus to show Theorem \ref{main} it suffices to show that $\chi [\Q [\Delta]/ \Theta,1]$ is a permutation character of $G$. A stronger result may be obtained when the order of $G$ is a prime power.

\begin{theorem} \label{mainprime} Let $\sigma$ and $G$ be as in Theorem \ref{main}, and let $\Delta$ and  $\Theta$ be as above. Then if the order of $G$ is a prime power, $\chi [\Q [\Delta]/ \Theta,q] = \sum_i \gamma_i q^i$, where each $\gamma_i$ is a permutation character of $G$.
\end{theorem}

\section{Permutation Characters of Cyclic Groups} \label{permchars}
 We begin by establishing a criterion for determining whether a character of a cyclic group is a permutation character. For basic character theory definitions and results see \cite{MR2270898}. 

Let $\chi$ be a rational character of finite cyclic group $G$. Now, $\chi$ is a permutation character of $G$ if and only if $\chi$ may be written as a linear combination of transitive permutation characters of $G$ with non-negative integer coefficients. Recall that the transitive permutation characters of $G$ are of the form $\ind{H}{G}$ for $H$ a subgroup of $G$, and note that as $G$ is Abelian

\begin{equation} \label{indchar}
\ind{H}{G}(g) = 
   \begin{cases}     [G:H] & \text{if } g \in H\\
      0 & \text{otherwise}.
   \end{cases}
\end{equation}

Now, as $\chi$ is rational, if $g$ and $g'$ in $G$ generate the same subgroup, then $\chi (g) = \chi(g')$. Hence for $H \leq G$ we may define $\chi (H) = \chi (h)$ where $h$ is any generator of $H$. Let $\mu$ be the M\"{o}bius function of the subgroup lattice $S(G)$ of $G$ (see, e.g., \cite{MR1442260}). Define $F_{\chi} : S(G) \rightarrow \R$ by
\begin{equation*}
 F_{\chi}(H) = \sum_{H \leq K} \mu(H,K) \chi(K).
\end{equation*}

\begin{lemma} \label{cexp} $\chi = \sum_{H \leq G} \frac{F_{\chi}(H)}{[G:H]}\ind{H}{G}$.
  \begin{proof}
    This is an immediate consequence of  (\ref{indchar}) and the M\"{o}bius inversion formula. For $g \in G$
    \begin{eqnarray*}
       \left( \sum_{H \leq G} \frac{F_{\chi}(G)}{[G:H]}\ind{H}{G} \right)(g) & = & \sum_{g \in H \leq G} \frac{F_{\chi}(G)}{[G:H]}[G:H]\\
        & = & \sum_{\gen{g} \leq H} F_{\chi}(H)\\
        & = & \chi(\gen{g})\\
        & = & \chi(g).
    \end{eqnarray*}

  \end{proof}

\end{lemma}
 Note that the M\"{o}bius inversion formula also implies that this is the unique way to write $\chi$ as a linear combination of transitive permutation characters. Thus $\chi$ will be a permutation representation of $G$ if and only if $\frac{F_{\chi}(H)}{[G:H]}$ is a non-negative integer for each $H \leq G$. It turns out that much of this condition will hold automatically for any rational character.

\begin{lemma} \label{cychar} Suppose $\chi$ is any rational character of $G = \Z_n$, and there exists $\alpha_H \in \mathbb{R}$ for each $H \leq G$ such that $\chi = \sum_{H \leq G}\alpha_H\ind{H}{G}$. Then for all subgroups $H$ of $G$, $\alpha_H \in \Z$. Furthermore, $\alpha_{\gen{1}} \geq 0$.

\begin{proof} Let $c$ be a generator of $G$, and let $\omega$ be a primitive $n^{th}$ root of 1. Let $\chi_i$ be the irreducible character of $G$ defined by $\chi_i(c) = \omega^i$. (Note that all of the irreducible characters of $G$ take this form for some $i \in \{ 1, \dots , n \}$.)
Any subgroup of $G$ may be expressed as $\gen{c^k}$, where $k$ is some divisor of $n$;  write $\gamma_k = \ind{\gen{c^k}}{G}$, $\alpha_k = \alpha_{\gen{c^k}}$. Now, $\chi$ has a unique decomposition into irreducible characters, $\chi = \sum_{i=1}^n \beta_i \chi_i$, where each $\beta_i$ is a non-negative integer. Observe that if $n | ik$, we have
\begin{eqnarray*}
 \ip{\gamma_k}{\chi_i} & = & \frac{1}{n}\sum_{g \in G} \gamma_k(g) \overline{\chi_i}(g)\\
  & = & \frac{1}{n}\sum_{g \in \gen{c^k}} k \overline{\chi_i}(g)\\
  & = & \frac{1}{n}\sum_{j = 1}^{\frac{n}{k}} k \omega^{-ijk}\\
  & = & \frac{1}{n}\sum_{j = 1}^{\frac{n}{k}} k \\
  & = & 1.
\end{eqnarray*}
On the other hand, if $n$ does not divide $ik$, $\omega^{-ik}$ is some primitive $l^{th}$ root of $1$, where $l > 1$ divides $\frac{n}{k}$. Then $\sum_{r = 1}^{l} \omega^{-ikr} = 0$, and so
\begin{eqnarray*}
 \ip{\gamma_k}{\chi_i} & = & \frac{1}{n}\sum_{j = 1}^{\frac{n}{k}} k \omega^{-ijk}\\
  & = & \frac{1}{n} \frac{n}{kl} k \sum_{r = 1}^{l} \omega^{-ikr} \\
  & = & 0.
\end{eqnarray*}
Now, for fixed $i | n $, $n$ divides $ik$ for a given divisor $k$ of $n$ if and only if $k = l \frac{n}{i}$ for some integer $l$. Furthermore, we then have $i = l \frac{n}{k}$, so $l | i$. Thus we see that $\beta_i = \sum_{l |i} \alpha_{l\frac{n}{i}}$.

Then $\beta_1 = \sum_{l |1} \alpha_{ln} = \alpha_n$. Thus, $\alpha_n$, the coefficient of $\ind{\gen{1}}{g}$, is a non-negative integer, proving the second assertion of the theorem.
Now we proceed inductively. Suppose $k|n$, and for all divisors $m$ of $n$ such that $k | m$, $k \neq m$, we have $\alpha_m \in \Z$. Let $j = \frac{n}{k}$ (so $j |n$). Then
\begin{eqnarray*}
\beta_j & = & \sum_{l|j} \alpha_{l\frac{n}{j}}\\
& = & \sum_{l|j} \alpha_{lk}\\
& = & \alpha_k + \sum_{l|j, l \neq 1} \alpha_{lk}
\end{eqnarray*}
Now, note that if $l \neq 1$ and $l |j$, then $k | lk$, $k \neq lk$, and $lk | n$, and so by induction $\sum_{l|j, l \neq 1} \alpha_{lk} \in \Z$, as is $\beta_j$, and thus $\alpha_k \in \Z$.
\\
\end{proof}
\end{lemma}

\section{Characterizations of $\chi [\Q [\Delta]/ \Theta,q]$}

 Throughout this section and the next, let $\sigma$ be a complete simplicial fan in $d$-dimensional real Euclidean space $V$ with associated simplicial complex $\Delta$, and let $G$ be a finite cyclic subgroup of $GL(V)$ which acts properly on $\sigma$ and acts as homomorphisms on the underlying lattice, $L$. Let $n=|G|$ and write $\chi_u = \chi [\Q [\Delta]/ \Theta,1]$.
 Our principal tool will be a formula for $\chi [\Q [\Delta]/ \Theta,q]$ due to Stembridge. For $g \in G$, let $V_g = \{ v \in V : gv = v \}$. Because $G$ acts properly on $\sigma$, the restriction of $\sigma$ to $V_g$ gives a complete simplicial fan $\sigma^g$ in $V_g$, and the simplicial complex associated to $\sigma^g$ is $\Delta^g = \{ F \in \Delta : gF = F \}$. Let $\delta (g) = \text{dim} (V_w)$, and for $\Gamma$ a simplicial complex, let $P_{\Gamma}(q)$ denote the h-polynomial of $\Gamma$.
 
\begin{lemma} \emph{\cite{MR1279220}} For $\Delta$, $G$ as above, and $g\in G$, $\chi [\Q [\Delta],q](g) = \frac{P_{\Delta^g}(q)}{(1-q)^{\delta (g)}}$.
\end{lemma}

We next need to describe our special system of parameters $\Theta$. Let $\varepsilon_1, \ldots , \varepsilon_n$ be a basis for $L$, and let $\ip{\cdot}{\cdot}$ be the standard inner product on $V$. For $v \in I$, let $x_v$ be the corresponding variable in $\Q [\Delta ]$. We then define
\begin{equation*}
\theta_i = \sum_{v \in I} \ip{v}{\varepsilon_i}x_v,
\end{equation*}
and let $\Theta = ( \theta_1, \ldots , \theta_n)$. Observe that for $g \in G$, we have
\begin{eqnarray*}
g\theta_i & = & \sum_{v \in I} \ip{v}{\varepsilon_i}x_{gv}\\
          & = & \sum_{v \in I} \ip{g^{-1}v}{\varepsilon_i}x_v\\
          & = & \sum_{v \in I} \ip{v}{(g^{-1})^*\varepsilon_i}x_v. 
\end{eqnarray*}
 Thus the $\Q$-span of $\theta_1, \ldots , \theta_n$ is $G$-stable, and in particular the representation of $G$ on this span and the representation of $G$ on $V$ are contragradient. But then as our vector space is real, they are in fact equivalent. In particular, $\Q [\Theta ]$ is isomorphic as a $\Q G$-module to the symmetric algebra of $V$, and thus 
\begin{equation*}\label{eqnparspan}
\chi [\Q [\Theta],q](g) = \text{det}(1-qg).
\end{equation*} 
 
 Finally, Stembridge showed in \cite{MR1279220} that $\Q [\Delta]$ and $\Q [\Delta] / \Theta \otimes \Q [\Theta]$ are isomorphic as $\Q G$-modules. It then follows that 
 \begin{equation*} \label{gchar}
\chi [\Q [\Delta]/ \Theta,q](g) = P_{\Delta^g}(q) \frac{\text{det}(1-qg) } {(1-q)^{\delta (g) }}.
\end{equation*}
In particular,
\begin{equation*} \label{uchar}
  \chi_u(g) = | \Delta^g | \text{det}_{V_g^{\bot}}(1-g),
\end{equation*}
 where $|\Delta^g |$ denotes the number of facets of $\Delta^g$.
 
For $g \in G$, write $Q_g(q) = \frac{\text{det}(1-qg) } {(1-q)^{\delta (g) }}$.

\begin{lemma} \label{qpoly} Suppose $g \in G$. Then for each divisor $l$ of $|g|$, there exists a non-negative integer $r_l$ such that every primitive $l^{th}$ root of unity appears as an eigenvalue of $g$ with multiplicity $r_l$, and
\begin{equation*}
Q_g(q) = \prod_{1\neq l | |g|}\Phi_l^{r_l}(q),
\end{equation*}
where $\Phi_l$ denotes the $l^{th}$ cyclotomic polynomial.
\begin{proof}
  Note that $\delta (g)$ is equal to the multiplicity of $1$ as an eigenvalue of $g$, and thus $Q_g(q) = \prod_i (1-\lambda_iq)$, where $( \lambda_i )_i$ is the list of eigenvalues (repeated according to multiplicity) of $g$ not equal to 1. Furthermore, as $g$ preserves the lattice $L$, the matrix for $g$ with respect to a basis for $L$ has rational coefficients, and hence $Q_g(q)$ has rational coefficients. Thus the complex eigenvalues of $g$ must come in conjugate pairs. The eigenvalues of $g$ are roots of 1, and thus each $\lambda_i$ is either $-1$ or appears with its conjugate, $\lambda_i^{-1}$. In particular it follows that $Q_g(q)$ is a monic polynomial with rational coefficients whose roots are exactly the eigenvalues of $g$ not equal to 1 (with corresponding multiplicities). The lemma then follows from the irreducibility of the cyclotomic polynomials over the rationals.
\end{proof}
\end{lemma}

In the next lemma we collect several properties of the cyclotomic polynomials which will be of use and are simple to check.

\begin{lemma} \label{cyclotomic}
Let $p$ be prime, $l$ a positive integer. 
\begin{enumerate} \item $\Phi_{p^l}(1) = p$. 
                  \item The degree of $\Phi_{p^l}$ is $p^{l-1}(p-1)$.
                  \item If $k > 1$ is not prime, then $\Phi_k(1) = 1$.
\end{enumerate}
\end{lemma}

Now let $c$ be a generator for $G$ and let $l | n$. The polynomial $Q_{c^l}(x)$ can be determined from $Q_c(x)$. By Lemma \ref{qpoly} we may write
\begin{equation*}
Q_c(q) = \prod_{k | n, k \neq 1} \Phi_k^{a_k}(q)
\end{equation*}
for some non-negative integers $a_k$. Now, the eigenvalues of $c^l$ are exactly the $l^{th}$ powers of the eigenvalues of $c$ (with corresponding multiplicities). If $\lambda$ is a primitive $k^{th}$ root of $1$, then $\lambda^l$ is a primitive $\frac{k}{(k,l)}^{th}$ root of 1. Thus
\begin{equation*}
Q_{c^l}(q) = \prod_{k|n, k \nmid l} \Phi_{\frac{k}{(k,l)}}^{a_kd(k,l)}(q),
\end{equation*}
where $d(k,l) = \frac{\text{deg}(\Phi_k)}{\text{deg}(\Phi_{\frac{k}{(k,l)}})}$.

In particular, let $M(p,N) = \text{max}\{ i : p^i | N \}$ for any prime $p$ and positive integer $N$. Then set $b(p,l,n) = 0$ if $M(p,l) = M(p,n)$, and otherwise 
\begin{equation*}
b(p,l,n) = \sum_{i=M(p,l)+1}^{M(p,n)}\sum_{s | \frac{l}{p^{M(p,l)}}} a_{p^is}d(p^is,p^{i-M(p,l)}).
\end{equation*}
Then by Lemma \ref{cyclotomic}
\begin{equation*}
Q_{c^l}(1) = \prod_{ p\text{ prime}}p^{b(p,l,n)}.
\end{equation*}

Note then that $\chi_u$ takes positive integer values, as each $b(p,l,n)$ is a non-negative integer. For $G$ having prime order, this is enough to imply that $\chi_u$ is a permutation character (this follows from the proof of Lemma \ref{cychar}, for example). To show Theorem \ref{main} we will need to restrict our attention to certain quotients of $G$, for which we will require the following lemma.

\begin{lemma} Suppose $l|n$. Then $G/\gen{c^l}$ acts properly on $\Delta^{c^l}$, as a subgroup of $GL(V_{c^l})$ and a set of homomorphisms on the restriction of $L$ to $V_{c^l}$, by letting $\overline{w}v = wv$ whenever $w$ is a class representative of $\overline{w} \in G/\gen{c^l}$. Furthermore, this gives a (graded) character, $\chi^l[q]$, for $G/\gen{c^l}$,  with
\begin{equation*}
\chi^l[q](\overline{c^j}) = P_{\Delta^{c^j}}(q)\prod_{k|l, k \nmid j} \Phi_{\frac{k}{(k,j)}}^{a_kd(k,j)}(q)
\end{equation*}
for each $j | l$. In particular,
\begin{equation*}
\chi_u^l := \chi^l[1](\overline{c^j}) = | \Delta^{c^j} | \prod_{p\text{ prime}}p^{b(p,j,l)}.
\end{equation*}

\begin{proof} That the action is well defined follows immediately from the fact that $\gen{c^l}$ acts trivially on $V_{c^l}$; the rest of the first claim then follows immediately from the properties of the original action.

Letting $\chi^l[q]$ be the character for the induced action on $\Q[\Delta^{c^l}]/ \Theta_l$ (where $\Theta_l$ is our special system of parameters for $V_{c^l}$ and $\Delta^{c^l}$), Theorem \ref{gchar} gives

\begin{equation*}
\chi^l[q](\overline{c^j}) = P_{(\Delta^{c^l})^{\overline{c^j}}}(q) \frac{\text{\emph{det}}_{V_{c^l}}(1-q\overline{c^j}) } {(1-q)^{\delta \overline{(c^j)} }}.
\end{equation*}
Furthermore, observe that if $j |l$, $\Delta^{c^j} \subset \Delta^{c^l}$, so $(\Delta^{c^l})^{\overline{c^j}} = \Delta^{c^j}$. Finally, it is clear that an eigenvalue of $c$ occurs as an eigenvalue of $\overline{c}$ if and only if its order divides $l$. Thus our earlier arguments, restricted now to $V_{c^l}$, yield the desired formulas.

\end{proof}

\end{lemma}

\begin{corollary} \label{eqrelate} For $j | l$ 
\begin{equation*} 
\chi [\Q [\Delta]/ \Theta,q](c^j) = \chi^l[q](\overline{c^j}) \prod_{k|n, k \nmid l, k \nmid j} \Phi_{\frac{k}{(k,j)}}^{a_kd(k,j)}(q).
\end{equation*}
Furthermore,
\begin{equation*} \label{relate}
\chi_u(c^j) = \chi_u^l(\overline{c^j})\prod_{p\text{ prime}}p^{b(p,j,n)-b(p,j,l)}.
\end{equation*}
\end{corollary}

We will require a slightly different characterization of the relationship between $\chi_u$ and $\chi_u^l$. To that end, suppose $p$ is prime, $l | n$, and $j |l$. Set $M = M(p,l)$ and $r = M(p,j)$, so $j = p^rk$ for some $k$ not divisible by $p$. If $M = M(p,n)$, set $b^l(p,j,n) = 0$. Otherwise, set 
\begin{equation*}
b^l(p,j,n)= \begin{cases} \sum_{i=M+1}^{M(p,n)}a_{p^ij}d(p^ij, p^i), &\text{if } r=0\\
                         \sum_{i=M+1}^{M(p,n)}(p-1)a_{p^ik}d(p^ik,p^{i-r-1}), &\text{otherwise}
           \end{cases}.
\end{equation*}
Finally, set
\begin{equation*}
C^l(a) = \prod_{p \text{ prime}} p^{b^l(p,a,n)}.
\end{equation*}
\begin{lemma} Let $j$, $l$ and $C^l$ be as defined above. Then
\begin{equation*}
\chi_u(c^j) = \chi_u^l(\overline{c^j})\prod_{a | j} C^l(a).
\end{equation*}
\begin{proof} This is equivalent to showing that
\begin{equation*}
\prod_{p\text{ prime}}p^{b(p,j,n)-b(p,j,l)} = \prod_{a | j} C^{l}(a),
\end{equation*}
i.e., for each prime $p$, $b(p,j,n) - b(p,j,l) = \sum_{a | j} b^l(p,a,n)$.
Note that if $p \nmid \frac{n}{j}$, $b(p,j,n) = b(p,j,l) = b^l(p,a,n) = 0$ (for all $a$). Thus it will suffice to show that $\sum_{a | j} b^l(p,a,l) = b(p,j,n)-b(p,j,l)$ for each $ p | \frac{n}{j}$.

By Lemma \ref{cyclotomic}, for $i > 0$, $\text{deg}(\Phi_{p^{i+1}}) = p \text{deg}(\Phi_{p^{i}})$. Thus, if $p^i | N$,
\begin{eqnarray*}
d(N,p^i) & = & \frac{\text{deg}(\Phi_N)}{\text{deg}(\Phi_{p^i})}\\
         & = & p\frac{\text{deg}(\Phi_N)}{\text{deg}(\Phi_{p^{i+1}})}\\
         & = & pd(N,p^{i+1}).
\end{eqnarray*} 
Furthermore, as $j | l$, $ r = M(p,j) \leq M$, so
\begin{eqnarray*}
b(p,j,n)-b(p,j,l) & = & \sum_{i=r+1}^{M(p,n)}\sum_{s | k} a_{p^is}d(p^is,p^{i-r)}) - \sum_{i=r+1}^{M}\sum_{s | k} a_{p^is}d(p^is,p^{i-r})\\
                  & = & \sum_{i=M+1}^{M(p,n)}\sum_{s | k} a_{p^is}d(p^is,p^{i-r}).
\end{eqnarray*}
But
\begin{eqnarray*}
\sum_{a | j}b^l(p,a,n) & = & \sum_{q = 0}^r \sum_{s | k} b^{l}(p,p^qs,n)\\
                       & = & \sum_{s | k}\left[ b^l(p,s,n) + \sum_{q=1}^rb^l(p,p^qs,n)\right] \\
                       & = & \sum_{s | k}\left[ \sum_{i=M+1}^{M(p,n)}a_{p^is}d(p^is, p^i)\right. \\
                       & & \left. + \sum_{q=1}^r\sum_{i=M+1}^{M(p,n)}(p-1)a_{p^is}d(p^is,p^{i-q-1})\right] \\
                       & = & \sum_{s | k}\sum_{i = M+1}^{M(p,n)}a_{p^is}\left[ d(p^is,p^i) + \sum_{q=1}^r(p-1)d(p^is,p^{i-q-1})\right] \\
                       & = & \sum_{s | k}\sum_{i = M+1}^{M(p,n)}a_{p^is}\left[ d(p^is,p^i) \right. \\
                       & & \left. + \sum_{q=2}^{r+1}d(p^is,p^{i-q-1}) - \sum_{q=1}^{r}d(p^is,p^{i-q-1}) \right] \\
                       & = & \sum_{s | k}\sum_{i = M+1}^{M(p,n)}a_{p^is}\left[ d(p^is,p^i) + d(p^is,p^{i-r}) - d(p^is,p^i)\right] \\
                       & = & \sum_{s | k}\sum_{i = M+1}^{M(p,n)}a_{p^is}d(p^is,p^{i-r})\\
                       & = & b(p,j,n) - b(p,j,l),
\end{eqnarray*}
as desired.
\end{proof}
\end{lemma}

\section{Proof of the Main Theorems}
We are now ready to prove Theorem \ref{main}. In light of Lemma \ref{cychar}, it will suffice to show that $F_{\chi_u}(H) \geq 0$ for all non-trivial subgroups $H$ of $G$. But each subgroup of $G$ is of the form $\gen{c^l}$ for some divisor $l$ of $n$. Thus we need only show that
\begin{eqnarray*}
0 & \leq & F_{\chi_u}(\gen{c^l})\\
  & = & \sum_{\gen{c^l} \leq H} \mu \left( \gen{c^l},H \right) \chi_u(H)\\
  & = & \sum_{j | l} \mu(\gen{c^l}, \gen{c^j}) \chi_u(c^j)\\
  & = & \sum_{j |l } \mu(\gen{c^l}, \gen{c^j})  \chi^l(\overline{c^j})\prod_{a | j} C^l(a),
\end{eqnarray*}
for every proper divisor $l$ of $n$, where $\chi^l$ and $C^l$ are defined as in the previous section. By Lemma \ref{cychar} we know that $F_{\chi_u^l}(\gen{\overline{1} }) = F_{\chi_u^l}\left( \gen{\overline{c^l} } \right)$ is non-negative, i.e.,
\begin{equation*}
\sum_{j |l } \mu_l(\gen{\overline{c^l}}, \gen{\overline{c^j}})  \chi_u^l(\overline{c^j}) \geq 0.
\end{equation*}
Here $\mu_l$ is the M\"{o}bius function on the subgroup lattice of $G / \gen{g^l}$. But this lattice is isomorphic to the interval $[\gen{c^l}, G]$ in $S(G)$, so
\begin{equation*}
\sum_{j |l } \mu_l(\gen{c^l}, \gen{c^j})  \chi_u^l(\overline{c^j}) \geq 0.
\end{equation*}

Furthermore, recall that in the case that the order of $G$ is prime, the non-negativity of $\chi$ immediately implies that $\chi$ is a permutation character. Thus we may induct on the number of factors of $n$, and assume that the restriction of $\chi_u^l$ to a proper subgroup of $G / \gen{g^l}$ gives a permutation character. In particular, for any fixed divisor $k$ of $l$, 
\begin{equation*}
\sum_{k | j, j |l } \mu_l(\gen{c^l}, \gen{c^j})  \chi_u^l(\overline{c^j}) \geq 0.
\end{equation*}

The proof of Theorem \ref{main} is thus completed by the following lemma:

\begin{lemma} Let $P$ be a finite bounded poset and suppose $f': P \rightarrow \R$ such that for all $p \in P$, $\sum_{p \leq q}f'(q) \geq 0$. Let $C$ be any map from $P$ into the positive integers, and define $f: P \rightarrow \R$ by $f(p) = f'(p)\prod_{q \leq p} C(q)$. Then $\sum_{p \in P}f(p) \geq 0$.

\begin{proof} We will first induct on $|P|$. If $|P| = 1$, the result is immediate. Suppose $|P| > 1$ and that the lemma holds for any bounded poset with size strictly less than $|P|$. Let $\hat{1}$ and $\hat{0}$ denote the greatest and least elements of $P$, respectively.

Extend the order on $P$ to a total order $\dot{\leq}$, and for $q \in P$, define $f_q(p) = f'(p)\prod_{q \dot{\leq} l \leq p}C(l)$ (so in particular $f = f_{\hat{0}}$). We will induct ``down'' with respect to $\dot{\leq}$ to show the result holds for each $f_q$.

By hypothesis $f'(\hat{1}) \geq 0$, so
\begin{eqnarray*}
\sum_{p\in P} f_{\hat{1}}(p) & = & f'(\hat{1})C(\hat{1}) + \sum_{p < \hat{1}}f'(p)\\
                             & \geq & \sum_{p\in P} f'(p)\\
                             & \geq & 0. 
\end{eqnarray*}
Now suppose $\hat{0} < q < \hat{1}$ and that for each $k \dot{\geq} q$, $\sum_{p\in P} f_{k}(p) \geq 0$. Let $r$ be the immediate successor to $q$ with respect to $\dot{\leq}$, and let $P' = \{ p \in P$ : $p \geq q \}$. Then $|P'| < |P|$ and $\sum_{p \leq l, l \in P'}f'(l) = \sum_{p \leq l, l \in P}f'(l) \geq 0$. Hence by our first induction hypothesis, $\sum_{q \leq p}f_r(p) \geq 0$. Furthermore, by our second induction hypothesis, $\sum_{p \in P}f_r(p) \geq 0$. Now, if $p \ngeq q$,
\begin{eqnarray*}
f_q(p) & = & f'(p)\prod_{q \dot{\leq} l \leq p}C(l)\\
       & = & f'(p)\prod_{r \dot{\leq} l \leq p}C(l)\\
       & = & f_r(p),
\end{eqnarray*}
while if $p \geq q$,
\begin{eqnarray*}
f_q(p) & = & f'(p)\prod_{q \dot{\leq} l \leq p}C(l)\\
       & = & f'(p)C(q)\prod_{r \dot{\leq} l \leq p}C(l)\\
       & = & f_r(p)C(q).
\end{eqnarray*}
Thus
\begin{eqnarray*}
\sum_{p \in P}f_q(p) & = & \sum_{p \geq q}f_q(p) + \sum_{p \ngeq q}f_q(p)\\
                     & = & C(q)\sum_{p \geq q}f_r(p) + \sum_{p \ngeq q}f_r(p)\\
                     & \geq & \sum_{p \geq q}f_r(p) + \sum_{p \ngeq q}f_r(p)\\
                     & = & \sum_{p \in P}f_r(p)\\
                     & \geq & 0.
\end{eqnarray*}
So by induction $\sum_{p \in P}f_q(p) \geq 0$ for all $q > \hat{0}$.

Finally, let $r$ now be the immediate successor to $\hat{0}$ in $\dot{\leq}$. Then
\begin{eqnarray*}
\sum_{p \in P}f_{\hat{0}}(p) & = &C(\hat{0})\sum_{p \in P}f_r(p)\\
                             & \geq & \sum_{p \in P}f_r(p)\\
                             & \geq & 0,
\end{eqnarray*}
 as desired.

\end{proof}

\end{lemma}

We now turn our attention to the case $|G| = p^r$ for some prime $p$ and positive integer $r$.  Here our computation is simplified as the subgroup lattice of $G$ is simply the total ordering $G = \gen{c} > \gen{c^{p}} >  \gen{c^{p^2}} \ldots  > \gen{c^{p^r}} = \gen{1}$. Then, for $j \leq i$,
\begin{equation*}
\mu \left( \gen{c^{p^i}}, \gen{c^{p^j}} \right) = \begin{cases} 1 & i=j\\
                                                   -1 & i = j+1\\
                                                   0 & \text{otherwise} 
                                    \end{cases}.
\end{equation*}

For ease of notation let $\chi[q](i) = \chi [\Q [\Delta]/ \Theta,q](c^{p^i})$.  Then, by Lemma \ref{cychar}, to show Theorem \ref{mainprime} it will suffice to show that $\chi[q](i) - \chi[q](i-1)$ has non-negative coefficients for $1 \leq i < r$.
 
First recall that $\Phi_{p^s}(q) = \sum_{i=0}^{p-1}q^{ip^{s-1}}$, and in particular has non-negative coefficients. Furthermore, $\Delta^{c^{p^s}}$ is a triangulation of a sphere and thus is Cohen-Macaulay, so $P_{\Delta^{c^{p^s}}}(q)$ has non-negative integer coefficients. In particular, $\chi^{p^s}[q](g)$ has non-negative coefficients for any $g \in G / \gen{c^{p^s}}$. Also, note that $d(p^l,p^{l-s}) = p^s$. 

Let $1 \leq i < r$. Then by Corollary \ref{eqrelate}
\begin{eqnarray} \label{powerform}
\lefteqn{\chi[q](i) - \chi[q](i-1)  = }\\ 
& & \chi^{p^i}[q](\overline{c^{p^i}})\prod_{l = i+1}^r \Phi_{p^{l-i}}^{a_{p^l}p^i}(q) \nonumber \\
& \text{}& -  \chi^{p^i}[q](\overline{c^{p^{i-1}}})\prod_{l = i+1}^r \Phi_{p^{l-i+1}}^{a_{p^l}p^{i-1}}(q). \nonumber
\end{eqnarray}
By Lemma \ref{cychar}, $\chi^{p^i}[q](c^{p^i}) - \chi^{p^i}[q](c^{p^{i-1}})$ has non-negative coefficients (this being the polynomials whose coefficients are those of $\ind{\gen{1}}{G / \gen{c^{p^s}}}$ in each degree), while $\chi[q](i) - \chi[q](i-1)$ is obtained from $\chi^{p^i}[q](c^{p^i}) - \chi^{p^i}[q](c^{p^{i-1}})$ by repeated multiplications of the first term by $\Phi_{p^r}^{p^i}(q)$ and the same number of multiplications of the second term by $\Phi_{p^{r+1}}^{p^{i-1}}(q)$.

Now, \begin{eqnarray*}
\Phi_{p^r}^{p^i}(q) & = & \left[ \left( \sum_{i=0}^{p-1}q^{ip^{r-1}} \right)^p \right]^{p^{i-1}} \\
                & = & \left[ \sum_{i=0}^{p-1}q^{ip^{r}} + R'(q) \right]^{p^{i-1}}\\
                & = & \Phi_{p^{r+1}}^{p^{i-1}}(q) + R(q),
\end{eqnarray*}
where $R'(q)$ and $R(q)$ are polynomials with positive coefficients. In particular, if $P$ and $Q$ are polynomials with non-negative coefficients such that $P - Q$ has non-negative coefficients, then
\begin{eqnarray*}
P(q)\Phi_{p^r}^{p^i}(q) - Q(q)\Phi_{p^{r+1}}^{p^i-1}(q) & = & P(q)\left[ \Phi_{p^{r+1}}^{p^{i-1}}(q) +R(q) \right] - Q(q)\Phi_{p^{r+1}}^{p^{i-1}}(q)\\
          & = & (P(q) - Q(q))\Phi_{p^{r+1}}^{p^{i-1}}(q) + P(q)R(q) 
\end{eqnarray*}
and thus has positive coefficients. Theorem \ref{mainprime} then follows from (\ref{powerform}).

\section{Remarks}
Theorem \ref{mainprime} does not hold for cyclic groups in general, that is, the representation of a cyclic group on $\Q [\Delta]/ \Theta$ need not be a permutation representation with respect to the grading. This follows from work of Stembridge in \cite{MR1279220}, where he notes that elements of certain Weyl groups acting as above yield graded character values with negative coefficients. Thus by restricting to the subgroup generated by such an element we obtain a character of a cyclic group that cannot be a permutation character in each graded component. Indeed we see that the arguments in the proof of Theorem \ref{mainprime} rely heavily on the fact that the cyclotomic polynomials for prime powers have non-negative coefficients, which is not true of cyclotomic polynomials in general.

It does seem possible that Theorem \ref{main} may be extended to larger classes of Abelian groups. Here the situation is complicated by the fact that there are in general many ways to write a character of a non-cyclic group as a linear combination of transitive permutation characters. Indeed, one may repeat the arguments of Section \ref{permchars}, now assigning any value at all to $\chi(H)$ for $H$ a non-cyclic subgroup of $G$ and keeping $\chi(\gen{g}) = \chi(g)$, and again obtain via M\"{o}bius inversion coefficients expressing $\chi$ as a linear combination of permutation characters. The problem then becomes one of whether values of $\chi(H)$ can be chosen for the non-cyclic subgroups in such a way that the resulting coefficients are non-negative integers. In certain simple cases (e.g., $Z_p \oplus Z_p$ and $\Z_{p^2} \oplus \Z_p$ for $p$ prime) this can be done more or less by hand, but a more general approach would be of interest.

Finally, note that although our character formulas involve terms coming from the face structure of fixed subcomplexes of $\Delta$, no knowledge of this structure was necessary in the proofs of our results (other than that the complexes are Cohen-Macaulay). This suggests the possibility of arriving at properties of $\Delta$ via investigation of the characters. In fact, work of this this flavor (using more directly the isomorphism of $H^*(X_{\sigma}, {\mathbb Q} )$ and  $\Q [\Delta]/ \Theta$ in order to use additional properties of the toric variety) has been done to investigate simplicial polytopes with certain simple symmetries, see, e.g.,  \cite{MR1354672}, \cite{MR2013971}, and \cite{MR932113}.

\end{document}